\documentclass[10pt, notitlepage]{article}
\usepackage[utf8]{inputenc}
\usepackage[english]{babel}
\usepackage{amsfonts}
\usepackage{amssymb}
\usepackage{amsmath}
\usepackage{amsthm}
\usepackage[titletoc]{appendix}
\usepackage{tocloft}
\usepackage{CJK}
\usepackage{mathrsfs}
\usepackage[english]{babel}
\usepackage{chngcntr}
\usepackage{graphicx}
\usepackage[all]{xy}
\newcommand{\widesim}[2][5]{ \mathrel{\overset{#2}{\scalebox{#1}[1]{$\sim$}}}
}
\counterwithin{equation}{section}

\makeatletter\@addtoreset{chapter}{part}\makeatother

\newcommand{\xdownarrow}[1]{%
  {\left\downarrow\vbox to #1{}\right.\kern-\nulldelimiterspace}
}

\begin{document}

\title{Cone construction via real intersection theory }

 \author{B. Wang (\begin{CJK}{UTF8}{gbsn}
汪      镔)
\end{CJK}}

\date {}

\renewcommand{\thefootnote}{\fnsymbol{footnote}}

\maketitle
 
\begin{abstract}
We show that the  cone construction  extends the Lefschetz standard conjecture
  to  the coniveau filtration.
\end{abstract}

\newcommand{\hookuparrow}{\mathrel{\rotatebox[origin=c]{90}{$\hookrightarrow$}}}
\newcommand{\hookdownarrow}{\mathrel{\rotatebox[origin=c]{-90}{$\hookrightarrow$}}}
\newcommand\ddaaux{\rotatebox[origin=c]{-90}{\scalebox{0.70}{$\dashrightarrow$}}} 
\newcommand\dashdownarrow{\mathrel{\text{\ddaaux}}}

\newtheorem{thm}{Theorem}[section]

\newtheorem{ass}[thm]{\bf {Claim} }
\newtheorem{prop}[thm]{\bf {Property } }
\newtheorem{prodef}[thm]{\bf {Proposition and Definition } }
\newtheorem{construction}[thm]{\bf {Main construction } }
\newtheorem{assumption}[thm]{\bf {Assumption} }
\newtheorem{proposition}[thm]{\bf {Proposition} }
\newtheorem{theorem}[thm]{\bf {Theorem} }
\newtheorem{apd}[thm]{\bf {Algebraic Poincar\'e duality} }
\newtheorem{cond}[thm]{\bf {Condition} }

\newtheorem{ex}[thm]{\bf Example}
\newtheorem{corollary}[thm]{\bf Corollary}
\newtheorem{definition}[thm]{\bf Definition}
\newtheorem{lemma}[thm]{\bf Lemma}
\newtheorem{con}[thm]{C}
\newtheorem{conj}[thm]{\bf Conjecture}
\newtheorem{result}[thm]{\bf Result}

\bigskip

\bigskip

\bigskip

\maketitle

\bigskip

\tableofcontents

\begin{center}\section{Introduction}\end{center}

\subsection {Generalization of the Lefschetz standard conjecture}
The Lefschetz standard conjecture were proposed by Grothendieck ([3]) in formulating a solution to Weil's conjectures.
The conjecture addresses a smooth projective variety $X$ of dimension $n$ over 
an algebraically closed field  of arbitrary characteristic.  We step back to assume the ground field is $\mathbb C$, and
the cohomology is Betti cohomology with rational coefficients.  Thus we consider the associated projective manifold which is still
denoted by $X$. Let $u$ be the hyperplane section class
in $H^2(X;\mathbb Q)$.  For $0\leq h\leq n$, let $L^{h}$ denote the homomorphism on the cohomology
\begin{equation}\begin{array}{ccc}
L^{h}: \sum_{i=0}^{2n-2h}H^{i}(X;\mathbb Q) &\rightarrow  & \sum_{i=0}^{2n-2h}H^{i+2h}(X;\mathbb Q)\\
\alpha &\rightarrow & \alpha\cdot u^{h}.
\end{array}\end{equation}
\bigskip

\bigskip

The hard Lefschetz theorem says $L^{h}$ has a property in topology. Precisely  the restriction  to $H^{n-h}(X;\mathbb Q)$ is 
an isomorphism to $H^{n+h}(X;\mathbb Q)$. 
Grothendieck envisioned that $L^h$ also has a property in  algebraic geometry. More precisely he proposed the Lefschetz standard conjecture:

\bigskip

\begin{conj} \par
 Let $A^j(X)\subset H^{2j}(X;\mathbb Q) $ for any interger $j$ be the subspace spanned by algebraic cycles.
Then the restriction  $L^{q-p}_0$ of  $L^{q-p}$ to $A^p(X)$,
\begin{equation}\begin{array}{ccc}
L^{q-p}_0:  A^p(X) &\rightarrow  & A^q(X)\\
\alpha &\rightarrow & \alpha\cdot u^{q-p}.
\end{array}\end{equation}
for $p+q=n, q\geq n$
is  an isomorphism.
\end{conj}

\bigskip

The conjecture  has been claimed  to be correct ([7]).  
In this paper, we use the real intersection theory  to extend
the same proof to the coniveau filtration.  This furthers Grothendieck's vision. 
\bigskip

The coniveau filtration is the decreasing filtration on the cohomology  $H^j(X;\mathbb Q)$ over $\mathbb Q$:
\begin{equation}
H^j(X;\mathbb Q)\supset \cdots\supset N^i H^j(X)\supset \cdots\supset N^{[{j\over 2}]}H^j(X)
\end{equation}
where $N^i H^j(X)$ is the subspace   spanned by 

\begin{equation}
Ker\bigl\{H^{j}(X;\mathbb Q) \to H^{j}(X- W;\mathbb Q)\bigr\}
 \end{equation}
for all subvarieties  $W$ of codimension at least $i$.  The index $i$ is called the coniveau and $j-2i$  the level. 

\bigskip

\begin{theorem} (Main theorem) \quad\par Let $p, q, k$ be whole numbers satisfying 
$$p+q=n-k, p\leq q. $$
 When $L^{q-p}$ are restricted to the subgroups of the coniveau filtration, the restricted homomorphism \begin{equation}\begin{array}{ccc}
 L^{q-p}_k: N^pH^{2p+k}(X) &\rightarrow  & N^qH^{2q+k}(X)\\
\alpha &\rightarrow & \alpha\cdot u^{q-p}.
\end{array}\end{equation} is an isomorphism, 
where $k$ is the level of the coniveau filtration.

\end{theorem}

\bigskip

{\bf Remark} The special case for $k=0$ is Conjecture 1.1.
 
\bigskip

\subsection{ Outline of the proof }\par

 \bigskip

The logic path of the proof is  identical to that of [7].  But we'll make a modification 
along the path to address the currents in real intersection theory ( [5], [6]).  Because many parts of this paper are duplicates  of [7],
  [5] and  [6],  to avoid the redundancy, we will only cite the results and proofs from the references without proofs. 
\bigskip

We use angle bracket $\langle\bullet\rangle$ to denote the object in the cohomology.
\bigskip

For  the indexes $p, q, k$ in Theorem 1.2,  we construct  a linear map 
$$\langle Con_{h, p}\rangle:   H^{2q+k}(X;\mathbb Q)\rightarrow H^{2p+k}(X;\mathbb Q)$$ which will be proved to
\par

(1) respect the coniveau structure as it maps the subspace
$N^pH^{2p+k}(X)$ \par\hspace{1CC} onto  $N^qH^{2q+k}(X)$.

(2)  and be the topological  
 inverse  of the map (1.2), i.e.  \begin{equation}\begin{array}{c}\langle Con_{h, p}\rangle\circ L_k^{q-p}=identity\\
L_k^{q-p}\circ \langle Con_{h, p}\rangle=identity.\end{array}
\end{equation}
\par

 The topological homomorphism $\langle Con_{h, p}\rangle$ is reduced from an extrinsic  operator $Con_{h, p}$ which 
is  a homomrophsim on  currents  constructed  with external data.  It contains information that are both topological and algebraic, 
extrinsic and intrinsic.    

\bigskip

\subsubsection{Cone family of cycles of currents}

The operator $Con_{n, p}$ is constructed from a family of currents, which begins as a family of algebraic cycles. \bigskip

Let $\mathbb C^{n+2}$ be a linear space over $\mathbb C$ with a basis
\begin{equation}
\mathbf e_0, \cdots, \mathbf e_{n+1}
.\end{equation}
Let $h$ be a natural number $<n$. Consider two subspaces 
\begin{equation}\begin{array}{c}
\mathbb C^{n+2-h}=span(\mathbf e_0, \cdots, \mathbf e_{n+1-h}), \\
 \mathbb C^{h}=span(\mathbf e_{n+2-h}, \cdots, \mathbf e_{n+1}).
\end{array}\end{equation}
Then  \begin{equation}
\mathbb C^{n+2-h}\oplus \mathbb C^{h}=\mathbb C^{n+2}.
\end{equation} 

Next we consider a variation of $\mathbb C^{h}$. 
Let $\mathbb C\cup \{\propto\}\simeq \mathbf P^1$ be the parameter space of the variation, denoted by $\Upsilon$, where $\propto$ is the infinity point of $\mathbf P^1$.    
The variation is defined as  
 \begin{equation} \mathbb C^{h}_z=span ( z\mathbf e_{n+2-h}-\mathbf e_0, \mathbf e_{n+3-h}, \cdots, \mathbf e_{n+1}), \quad 
for \ z\in\mathbb C\end{equation}
and $\mathbb C^{h}_\propto$ is the original $\mathbb C^h$ which is 
the limit of subspaces $\mathbb C^{h}_z$ in Grassmannian  as $z\to\propto$.  
Let $U=\mathbb C^\ast\cup \{\propto\}$ be the affine open set that parametrizes those  $\mathbb C^{h}_z $ satisfying
\begin{equation}
\mathbb C^{n+2}=\mathbb C^{n+2-h}\oplus  \mathbb C^{h}_z.
\end{equation}
The only point $z=0$ not in $U$ corresponds to the plane $\mathbb C^{h}_0$ that fails the decomposition (1.11). We call $z=0$ the unsteady point, others steady points.  
Therefore for each steady point $z\in U$, we have the unique decomposition (1.11) 
\begin{equation}
\mathbf x=(\mathbf x_1(z), \mathbf x_2(z))
\end{equation}
for  the vector $\mathbf x\in \mathbb C^{n+2}$. 
  This decomposition  gives a regular map 
\begin{equation}\begin{array}{ccc}
 \mathbb C\times U\times (\mathbb C^{n+2-h}\oplus \mathbb C^{h}_z) &\rightarrow & 
\mathbb C^{n+2-h}\oplus  \mathbb C^{h}_z=\mathbb C^{n+2}\\
(t, z, (\mathbf x_1(z), \mathbf x_2(z))) &\rightarrow & (\mathbf x_1(z), t \mathbf x_2(z)).\end{array}\end{equation}
After completing  $\mathbb C,  U$ and projectivizing the linear spaces,  the regular map  determines a rational map of the projective variety
\begin{equation}\begin{array} {ccc}
\kappa: \mathbf P^1\times \Upsilon\times \mathbf P^{n+1} &\dashrightarrow & \mathbf P^{n+1}\\
(t, z, [\mathbf x_1(z), \mathbf x_2(z)]) &\dashrightarrow &[\mathbf x_1(z), t\mathbf x_2(z)],
\end{array}\end{equation}
where $t, z$ are points in the affine open sets $\mathbb C, U$.
Let \begin{equation}\Omega=graph(\kappa) \subset \mathbf P^1\times \Upsilon\times \mathbf P^{n+1}\times \mathbf P^{n+1}\end{equation}
where  the graph of a rational map is defined to be the closure of the  graph at the regular locus (   
the same  for images and preimages of rational maps).\bigskip

Now we consider  the smooth projective variety $X$ of dimension $n$. Let 
$$X\stackrel{\mu}\rightarrow \mathbf P^{n+1}$$ be a birational morphism to a hypersurface of $ \mathbf P^{n+1}$ with
very ample line bundle $\mu^\ast (\mathcal O_{\mathbf P^{n+1}}(1))$, and proper intersections with all the coordinates planes
 of the basis $\mathbf e_0, \cdots, \mathbf e_{n+1}$. Also $\mu(X)$ does not contain lines.  
The collection of above spaces $\mathbb C^{n+2}$, $\mathbb C^{n+2}, \mathbb C^h_z$ and $\mu$ is called {\bf cone data}. 
Next we construct a family of algebraic cycles with the cone data.  
Using $\Omega$ we define the subvariety $\Sigma$
to be the image of the rational map $(id, id, \mu^{-1}, \mu^{-1})=\tau^{-1}$,
\begin{equation}\begin{array}{ccc}
\Omega\cap (\mathbf P^1\times\Upsilon\times \mu(X)\times \mu(X))
&\dashrightarrow & \mathbf P^1\times\Upsilon\times X\times X,
\end{array}\end{equation}
where $\tau=(id, id, \mu, \mu)$ is the birational-to-image map.
We have a straightforward assertion for $\Sigma$.  
\bigskip

\begin{proposition}
$\Sigma$ is reduced with two components of the same dimension $n+1$, and one of components  
is
\begin{equation}
\{1\}\times \Upsilon\times \Delta_{X}
\end{equation}
where $\Delta_{X}$ is the diagonal of $X\times X$.
\end{proposition}

\bigskip

The proposition gives a definition
\bigskip

\begin{definition} \quad\par 
 We define  
\begin{equation}
\Theta
\end{equation}
to be the other component of $\Sigma$, i.e.
\begin{equation}
\Theta=\Sigma-(\{1\}\times \Upsilon\times \Delta_{X}).\end{equation}
  \end{definition}
\bigskip

We denote the algebraic cycle of $\Theta$ also by $\Theta$, and denote the fibres of $\Theta$ over  $t\in \mathbf P^1, z\in \Upsilon$
as cycles in $X\times X$ by $\Theta_t^z$. They are not always irreducible. Similarly $\Theta_t$  over each $t\in \mathbf P^1$ is a well-defined family of algebraic cycles in $\Upsilon\times X\times X$. Also to avoid the overwhelming notations in the context,  
the fibres lifted to  $\mathbf P^1\times \Upsilon\times X\times X$  are  denoted by
the same notations  $\Theta_t^z, \Theta_t$ respectively.  Their corresponding images in $\mathbf P^{n+1}\times \mathbf P^{n+1}$ and 
$\Upsilon\times \mathbf P^{n+1}\times \mathbf P^{n+1}$ are denoted by $\Omega_t^z$ and $\Omega_t$ respectively.   \bigskip

The intersection of the family members with an algebraic cycles led to the proof of Lefschetz standard conjecture in [7]. In this paper, we show that the same intersection with singular cycles yields Main theorem.  So next we use a novel tool that allows us to  take the intersection   with singular cycles.  
The tool is the real intersection theory that has been developed in [5], [6].  It is a study of 
geometry in algebraic  varieties  without the structure of a manifold through the geometric measure theory.
 The importance of real intersection theory 
goes beyond the content of this paper.  However we should 
give a brief introduction in the following.   On any manifold $Y$, 
there is a special type of currents inside of $\mathscr D'(Y)$, 
 called Lebesgue currents which  include  singular chains and $C^\infty$-forms, and they form a subspace denoted
by $\mathcal C(Y)$.  
   For any two Lebesgue currents  $T_1, T_2$,  we  define an extrinsic 
intersection current as the currents' limit of the De Rham's homotopy regularization ([1]), denoted by
\begin{equation}
[T_1\wedge T_2].
\end{equation}
The intersection is also Lebesgue,  but it depends on a special type of extrinsic covering 
 on the manifold called De Rham data.  
This currents' intersection,  even though is extrinsic, 
 but unites  all the known products.   
Real intersection theory is an extension of  Fulton's intersection theory ([2]).   For instance  the notion of a correspondence
is parallel to that in Fulton's intersection theory, and plays the key role.   Let $X\times Y$
be the Cartesian product of compact manifolds equipped with a De Rham data, and $\mathcal J\in \mathcal C(X\times Y)$ a Lebesgue current. Then there is a homomorphism $\mathcal J_\ast$ called the correspondence of   currents, defined by
\begin{equation}\begin{array}{ccc}
\mathcal J_\ast: \mathcal C(X) &\rightarrow & \mathcal C(Y)\\
\sigma &\rightarrow & (P_Y)_\ast [\mathcal J \wedge (\sigma\otimes Y)],
\end{array}\end{equation}
where $P_Y: X\times Y\to Y$ is the projection.  
In algebraic geometry, we assume $\mathcal J, \sigma$ are algebraic cycles and $\wedge$ is the well-defined
intersection of algebraic cycles. Then $\mathcal J_\ast$ is the usual correspondence. This consistency
provides the basis for our  study of  coniveau filtration.

\bigskip

{\bf Remark} Most of notions for real intersection theory resemble those in Fulton's intersection theory.  
However the dependence of  the extrinsic De Rham data  distinguish them
from each other.

\bigskip

Let's come back to the setting of the cone data.  First we let 
all spaces be equipped with De Rham data where the Cartesian product 
has product De Rham data, and each factor  has projection 
De Rham data (in order  to use the projection formula).   
Let $\sigma$ be a  singular cycle in $X$, which is also a Lebesgue current.  
Next we work with the real intersection theory.
For each $t\in \mathbf P^1$, we define $\Psi_t(\sigma)$ to be a family of currents
\begin{equation}
(\eta_4)_\ast [\Theta\wedge (\{t\}\otimes \Upsilon\otimes \sigma\otimes X)]\end{equation}
where $\eta_4: \mathbf P^1\times \Upsilon\times X\times X(4th\ factor)$ is the projection to the last factor $X$.
  Each member $\Psi_t(\sigma)$ is a closed current called $t$-end cycle. 
The main theorem is a result of careful analysis of three EXTRINSICALLY 
determined end cycles: $\Psi_0(\sigma)$, $\Psi_1(\sigma)$, $\Psi_\infty(\sigma)$. 
(However  in application they'll  not have  a common $\sigma$). \bigskip

 Real intersection theory shows
\bigskip

\begin{proposition}  All end cycles above are homotopic. \end{proposition}

\bigskip

This is because the calculation of intersection of currents for any two points 
$t_1, t_2\in \mathbf P^1$ leads to the homotopy
 \begin{equation}
\Psi_{t_1}(\sigma)-\Psi_{t_2}(\sigma)=d \Lambda
\end{equation}
where $\Lambda$ is a current of higher dimension and $d$ is the differential of currents. \bigskip

\noindent $\bullet$\quad {\bf 1-end cycle}\bigskip

\begin{proposition} Let $\mathcal C_0(X)$ be the subspace consisting of closed Lebesgue currents.  There is an algebraic cycle $\omega\in Z_n(X\times X)$ extrinsically determined  such that
\begin{equation}
\Psi_1(\sigma)=m\sigma+\omega_\ast (\sigma)\end{equation}
where $m$ is some natural number.
\end{proposition}
\bigskip

\noindent  $\bullet$ \quad {\bf $\infty$-end cycle}\bigskip

This end cycle when restricted to a Zariski open set gives the factorization of the identity in (1.6). 
It however does not have a straighforward process. \medskip

Let $$V^h=div (\mu^\ast(x_{n+1-h}))\cap\cdots\cap div (\mu^\ast(x_{n+1}))$$ be the  
$h$-codimensional,   smooth, irreducible plane section of
$X$ by the very ample line bundle $\mu^\ast (\mathcal O_{\mathbf P^{n+1}}(1))$.
The decomposition (1.11) has a natural projection, 
\begin{equation}\begin{array}{ccc}
 \mathbb C^{n+2-h}\oplus \mathbb C^h_z &\rightarrow &  \mathbb C^{n+2-h}, 
z\in U\end{array}\end{equation}
which gives the rational map, 
\begin{equation}\begin{array}{ccc}
\{\infty\}\times U\times \mathbf P^{n+1} &\dashrightarrow &  \mathbf P^{n+1-h}\\
(\infty, z, [\mathbf x]) &\rightarrow & [\mathbf x_1(z)]
\end{array}\end{equation}
$\mathbf x=\mathbf x_1(z)\oplus \mathbf x_2(z)$ is the unique decomposition (1.11).
Let $G$ be the transpose of its graph (transposed in the order $\{\infty\}\times U\times  \mathbf P^{n+1-h}\times  \mathbf P^{n+1}$). 
Let $\overset{\circ}I_h$ be the intersection cycle 
\begin{equation}
   (\{\infty\}\times U\times  V^h\times X)\cdot G.
\end{equation}

We denote the linear operator (for currents),
\begin{equation}\begin{array}{ccc}
\mathcal C_o(V^h) &\rightarrow &\mathcal C_o(X)\\
\alpha &\rightarrow & {1\over m+d}(\eta_4)_\ast\overline  {\biggl[\overset{\circ}{I_h}
\wedge (\{\infty\}\otimes U\otimes \alpha \otimes X)\biggr ] }
\end{array}\end{equation} by
$Con_h$,  where  the closure is defined as the convergence of currents (see Definition 3.6),  and 
 $d=deg(\mu^\ast(\mathcal O_{\mathbf P^{n+1}}))$ 

\bigskip

Let $v^h$ be the linear operator
$$\begin{array}{ccc}
\mathcal C_o(X) &\stackrel{v^h}\rightarrow & \mathcal C_o(V^h)\\
 \sigma &\rightarrow & [\sigma\wedge  V^h]_{V^h},\end{array}$$
 where $ [\sigma\wedge  V^h]_{V^h}$ is the intersection of currents of manifold $V^h$.

\begin{proposition}
Then the end cycle 
at infinity has a factorization, 
\begin{equation}
\Psi_\infty(\sigma)=(m+d)  Con_h\circ v^h (\sigma)+\mathcal T_\sigma
\end{equation}
for $\sigma\in \mathcal C_o(X)$ where  $\mathcal T_\sigma \in \mathcal C(X)$
is homologous to zero.  
\end{proposition}
\bigskip

{\bf Remark}
The factorization in (1.29)  has the direct link to the factorization  in (1.6).  However the formula holds for closed currents, not for 
the cohomology classes.

\bigskip

\noindent $\bullet$ \quad {\bf $0$-end cycle}

\bigskip

\begin{proposition} Let $\sigma\in \mathcal C_o(X)$. 
 If $dim(\sigma)<2n-2h$, then
\begin{equation}
\Psi_0(\sigma)
\end{equation}
is homologous to a current supported on $V^h$.
\end{proposition}

\bigskip

Let $H_{p, V^h}(X;\mathbb Q)$ be the homology supported on $V^h$. Then there is 
a sequence 
\begin{equation}\begin{array}{ccccc}
H_{p}(V^h;\mathbb Q) &\stackrel{\nu_1}\rightarrow & H_{p, V^h}(X;\mathbb Q)
&\stackrel{\nu_2}\rightarrow & H_{p}(X;\mathbb Q)\end{array}\end{equation}
with the surjective $\nu_1$. 
 Proposition 1.8    leads to a corollary, \bigskip

\begin{corollary} The map $\nu_2$ is surjective for $p<2n-2h$. Therefore 
the inclusion map which is the composition of (1.31), 
$$\begin{array}{ccc}
H_{p}(V^h;\mathbb Q) &\rightarrow & H_{p}(X;\mathbb Q)\end{array}
$$
 is also surjective.

\end{corollary}

\bigskip

{\bf Remark} The  corollary is a result in cohomology. It requires a cohomological descend, Proposition 1.10.   The corollary extends 
Lefschetz hyperplane theorem to smaller dimensions.   However the assertion is weaker on the larger dimensions. 
The proof of Corollary 1.8 does not resemble the currently known proofs of Lefschetz hyperplane theorem. 

\bigskip

\subsubsection{Two different aspects of end cycles}

\bigskip

End cycles have two contrasted aspects that yield the main theorem. First is the cohomological descend, which is purely topological. 
 The second is algebro-geometric  about the preservation of the level of coniveau filtration.  \bigskip

$\bullet$ \quad {\bf Cohomological descend}.  
First we consider the 1-end cycle. In transcendental geometry, we will have
\begin{equation}
\langle\omega_\ast (\sigma )\rangle\cup u= d \langle\sigma \rangle\cup u
\end{equation}
where $u$ is the hyperplane section class $c_1(\mu^{\ast}(\mathcal O_{\mathbf P^{n+1}}))$ and $d=deg(\mu^{\ast}(\mathcal O_{\mathbf P^{n+1}}))$. 
The hard Lefschetz theorem implies an equality in cohomology, 
 \begin{equation}
\langle\omega_\ast (\sigma )\rangle=d \langle\sigma \rangle
\end{equation}
for   $deg(\sigma)>{n}$. By the symmetry, (1.33) is extended to all non-middle dimensional cohomology groups.

\bigskip

\begin{proposition}
   So  for $\sigma\in \mathcal C_o(X)$ and $dim(\sigma)\neq {n}$, 
\begin{equation}\begin{array}{c}
{1\over m+d}\langle\Psi_1(\sigma)\rangle=\langle \sigma \rangle\end{array}\end{equation}
Therefore \begin{equation}
{1\over m+d}\langle\Psi_0(\sigma)\rangle=\langle\sigma\rangle.\end{equation}
\end{proposition}

\bigskip
{\bf Remark} In application, (1.34) and (1.35) may  not be applied to  a common $\sigma$.
\bigskip

 Now we study the descend of  $\infty$-end cycle. By the proposition 1.7, it amounts to prove the cohomological descend of
   $Con_h$. Nevertheless this descending in general is false due to the existence of primitive cycles in cohomology.   So we show the cohomological descend for a specific case. 
  We  restrict $Con_h$ to the closed Lebesgue currents 
of dimension ${n-h}$ on $V^h$.  We denote the restriction map by $Con_{h, p}$ where $p=n-h$.
Applying the real intersection theory, we obtain that 
\begin{equation}
v^h\circ Con_{h, p}(\delta)= \delta+\mathcal J_\delta,      
\end{equation} for $\delta\in \mathcal C_{o}(V^h)$, where $\mathcal J_\delta$ is a current homologous to $0$.   
The hard Lefschetz theorem on the equality (1.36) asserts
that $$\delta \widesim[8]{homological \ equi} 0 \ in \ X\Longrightarrow  Con_{h, p}(\delta) \widesim[8]{homological \ equi} 0\ in\ X.$$
Hence $Con_{h, p}$ descends to the cohomology as an isomorphism on the cohomology
\begin{equation}\begin{array}{ccc}
\langle Con_{h, p}\rangle:  H^{n+h}_{V^h}(X;\mathbb Q) &\simeq  & H^{n-h}(X;\mathbb Q).
\end{array}\end{equation}
 where 
$H^{n-h}_{V^h}(X;\mathbb Q)$ denotes the cohomology with the support on $V^h$. 
\bigskip

\noindent $\bullet$ \quad {Level of coniveau filtration}

It is obvious that $\langle Con_{h, p}\rangle$ represents the Poincar\' e
 duality  which does not  contain the algebro-geometric information. But the specific isomorphism $\langle Con_{h, p}\rangle$ does. 
In general the intersection of currents is additive on the index - level. Since the $Con_{h, p}$ is obtained by intersecting with the algebraic cycles which has level 0, $Con_{h, p}$ preserves the level in algebraic geometry. So does its cohomological descend.   More precisely using real intersection theory we obtain that $\langle Con_{h, p}\rangle$  sends 
the level $k$ cycles onto the level $k$ cycles for any $k$. 
So

\bigskip

\begin{proposition}
\begin{equation}\begin{array}{ccc}
\langle Con_{h, p}\rangle:  N^qH^{2q+k}_{V^h}(X;\mathbb Q) &\rightarrow   & N^pH^{2p+k}(X;\mathbb Q)

\end{array}\end{equation}
is an isomorphism where $p+q=n-k$.
\end{proposition}
\bigskip

This is the second aspect of end cycles -- algebraicity. 
\bigskip

At last we use the same homotopy for the end cycle at $0$, specifically Corollary 1.9,  to obtain  the isomorphism
\begin{equation}
 H^{2q+k}_{V^h}(X;\mathbb Q)\simeq H^{2q+k}(X;\mathbb Q)
\end{equation}
and furthermore the isomorphism preserves the level of cycles as before. 
Therefore there  is an isomorphism 

\begin{proposition}
 $$N^qH^{2q+k}_{V^h}(X;\mathbb Q)\simeq N^qH^{2q+k}(X;\mathbb Q).$$
\end{proposition}
\bigskip

So we obtain Main theorem  
\bigskip

\begin{equation}
N^qH^{2q+k}(X;\mathbb Q)\backsimeq N^pH^{2p+k}(X;\mathbb Q).
\end{equation}

\bigskip

\bigskip

In the following sections we give the details.  In section 2, we construct the cone family of currents. It starts 
in subsection 2.1, where we give a construction of a family of algebraic cycles. In subsection 2.2, we intersect the algebraic cycles with 
singular cycles to obtain a family of currents.  Main theorem comes from the study of three particular members of the family called end cycles. Thus in section 3, we study the  end cycles as currents.  In sections 4, we address two contrasted aspects of end cycles -- algebraicity and cohomologicity that result in Main theorem.

\bigskip

\bigskip

\bigskip

\begin{center}\section{Cone family}\end{center}

The cone family of currents is  a homotopy deduced from an algebro-geometric  deformation with parameter space $\mathbf P^1$. 
So first we introduce the algebraic deformation.
\bigskip

\subsection {Family of algebraic cycles}

We continue with the cone data established in section 1. In this subsection we prove Propositions 1.3 which is 
the understanding of the family of algebraic cycles. 
See [7] for the proof of Proposition 1.3.

\bigskip

\subsection{ Family of  cycles of currents}
\bigskip

\begin{proof} of Proposition 1.5: Let $\epsilon, 0\in \mathbf P^1$.   Let $S^1\subset \mathbf P^1$ be a real circle that pass
through $\epsilon, 0$.  Continuing from the setting of Proposition 1.5, we 
defined the current \begin{equation}
\mathcal J=[\Theta\wedge (S^1\otimes \Upsilon\otimes X\otimes X)]\end{equation}
which is closed. 
Now we use the diffeomorphism between a piece of $S^1$ and the closed interval $[0, \epsilon]$.
Applying Proposition 4.12, [6], we obtain that
\begin{align}\begin{split}
 &[ \mathcal J_\epsilon\wedge (\Upsilon\otimes \sigma\otimes X)]-[ \mathcal J_0\wedge  (\Upsilon\otimes \sigma\otimes X)]
\\& =d (-1)^{kp }((P_\mathcal X)_\ast \mathcal J_{I_\epsilon }(T))
\end{split}\end{align}
where $k=deg(\mathcal J), p=deg(T)$, $P_{X}:\mathbf P^1\times  X\to  X$ is the projection. 
Hence $$[ \mathcal J_\epsilon\wedge (\Upsilon\otimes \sigma\otimes X)]$$ is homotopic to 
$$[ \mathcal J_0\wedge  (\Upsilon\otimes \sigma\otimes X)].$$
By the associativity of currents intersection, Property 2.6, [6], we have
$$[ \mathcal J_\epsilon\wedge (\Upsilon\otimes \sigma\otimes X)]=[\Theta\wedge (\{\epsilon\}\otimes \Upsilon\otimes \sigma\otimes X)] $$ 
$$[ \mathcal J_0\wedge (\Upsilon\otimes \sigma\otimes X)]=[\Theta\wedge (\{0\}\otimes \Upsilon\otimes \sigma\otimes X)].$$
This completes the proof.\end{proof}

\begin{center}\section{End cycles in the cone family}\end{center}

\subsection{Algebraic part of  end cycles} To sufficiently understand the end cycles, we must have detailed analysis on the family of algebraic cycles $\Theta$. 
To do that, we use the following coordinates from cone data. 
Let $x_0, \cdots, x_{n+1}$ be the coefficients of the basis $\mathbf e_0, \cdots, \mathbf e_{n+1}$ for $\mathbb C^{n+1}$.
Then $x_0, \cdots, x_{n+1}$ are homogeneous coordinates for $\mathbf P^{n+1}$.  
Let  $z\neq 0$ be
complex numbers parametrize   the affine neighborhood of $\Upsilon$. 
 Then the homogeneous  coordinates  for $\mathbf P^{n+1-h}$, $\mathbf P^{h-1}_z$ as in the decomposition (1.11) are
\begin{equation}\begin{array}{c}
\mathbf x_1(z): x_0+{x_{n+2-h}\over z}, x_1, \cdots, x_{n+1-h}\\
\mathbf x_2(z):   {x_{n+2-h}\over z}, x_{n+3-h},  \cdots, x_{n+1}. \end{array}\end{equation}
In the product $$\mathbf P^1\times \Upsilon\times \mathbf P^{n+1}\times \mathbf P^{n+1},$$
the coordinates for the third factor $\mathbf P^{n+1}$ has $x$-coordinates as above.  
The same coordinates for the last factor $\mathbf P^{n+1}$ in the product
will be denoted by the letter $y$, 
 $$\begin{array}{c}
\mathbf y_1(z): y_0+{y_{n+2-h}\over z}, y_1, \cdots, y_{n+1-h}\\
\mathbf y_2(z):   {y_{n+2-h}\over z}, y_{n+3-h}, \cdots, y_{n+1}. \end{array}$$

 \smallskip

In  the following we focus on the component of its fibre at three points $t=0, 1, \infty$. 

\bigskip

\noindent $\bullet$ Case $t=1$.  \medskip

Let's  to find the algebraic cycle
$$\Theta^z_1\subset \{1\}\times U\times X\times X,$$
i.e. $z$ is steady. 
We push the cycles  forward to 
the projective spaces
$$\{1\}\times U\times \mathbf P^{n+1}\times \mathbf P^{n+1}$$
to intersect with the subvariety (of codimension 2),
$$\{1\}\times U\times \mu(X)\times \mu(X).$$

Notice $\mu(X)$ is a hypersurface of $\mathbf P^{n+1}$. Assume 
$\mu(X)$ is defined by a polynomial $f$.  Then $(\mu^2)_\ast (X\times X)$ is a complete intersection 
 defined by two polynomials $f(x), f(y)$ in $$\mathbf P^{n+1}\times \mathbf P^{n+1}.$$
where $x, y$ are coordinates for the two copies of $\mathbf P^{n+1}$ in the Cartesian product.

Then $\Omega_t^z\cap \mu^2(X\times X)$,  denoted by $\Theta_t^z$  is explicitly defined by
$$ f(x_1, t x_2)=f(x_1, x_2)$$
where $t\in \mathbb C$ and $x_i=x_i(z)$ as the coordinates split  in (3.1).   Observe the expansion 
$$f(x_1, tx_2)-f(x_1, x_2)=(t-1)^r g_r^z+(t-1)^{r+1} g_{r+2}^z+\cdots.$$
Then the specialization  as $t\to 1$ in $\Delta_{\mathbf P^{n+1}}$ of coordinates $[x_1, x_2]$ is 
defined  two polynomials $$f(x_1, x_2)=g_r^z(x_1, x_2)=0.$$
Therefore 
$\Theta_{1}^z$ is the specializtion which is equal to some hypersurfaces $\{g_r^z=0\}$ of the divisor 
$$\{f=0\}\simeq X\subset \Delta_{\mathbf P^{n+1}}$$ of the projective space.   
Thus the projection  of $\Theta_1|_{U\times X\times X}\subset U\times X\times X$ to $X\times X$ has the closure that equals to the diagonal 
 $\Delta_{X}$.  If $\mathcal P: \Upsilon\times X\times X\to X\times X$ is the projection, then as an algebraic cycle, 
\begin{equation}
\overline{\mathcal P_\ast (\Theta_1|_{U\times X\times X})}=m\Delta_{X},
\end{equation}
where $m$ is some natural number, and the closure is defined as the cycle projection of the graph of the rational map, 
$$\mathcal P|_{U\times X\times X}:   \Theta_1|_{U\times X\times X}\to X\times X.$$ 
\bigskip

Now we consider the unsteady point $\propto$. The fibre $\Omega_1^{\propto}$ is an irreducible subvariety of dimension $n+2$ in $\mathbf P^{n+1}\times \mathbf P^{n+1}$ defined by

\begin{equation} \begin{array}{c}
x_i y_j=x_j y_i, 1\leq i, j\leq n+1
\end{array}
\end{equation}
Let's denote it by $K_1$.
So $\omega=K_1\cdot (X\times X)$ has dimension $n$. This  is a part of the cycle $\Theta_1$, originated
from the fibre in $\{\propto\}\times \{1\}\times X\times X$.  Hence we have

\bigskip

\begin{proposition}
\begin{equation}
\mathcal P_\ast (\Theta_1)=m\Delta_{X\times X}+\omega.
\end{equation}
\end{proposition}
\bigskip

\noindent $\bullet$  Case $t=\infty$. \medskip

$\Omega_\infty$ is a subvariety defined by
\begin{equation}\left\{ \begin{array}{c}
x_i y_j=x_j y_i, 1\leq i, j\leq n+1-h\\
x_i y_j=x_j y_i, n-h+2\leq i, j\leq n+1\\
(zx_0+x_{n+2-h}) y_j=(zy_0+y_{n+2-h}) x_j, 1\leq j\leq n+1-h\\
 x_i y_j =0, 1\leq j\leq n+1-h, n-h+2\leq i\leq n+1, \\
x_j (zy_0+y_{n+2-h}) =0, n-h+2\leq j\leq n+1.
\end{array}\right.
\end{equation}  
\bigskip

First we consider the open component in
$$\{\infty\}\times U \times X\times X,$$
surjective to $U$, 
i.e, considering generic $z\neq 0$. 
Observing the fourth equations $x_i y_j =0$, we obtain that $\Omega_\infty$ over generic for $z\in U$ has two components:
\par
$K_2$ defined by
 \begin{equation}\left\{ \begin{array}{c}
x_i=0, n-h+2\leq i\leq n+1,\\
x_i y_j-x_j y_i=0, 1\leq i, j\leq n+1-h\\
x_i y_j-x_j y_i=0, n-h+2\leq i, j\leq n+1\\
(zx_0+x_{n+2-h}) y_j-(zy_0+y_{n+2-h}) x_j=0, 1\leq j\leq n+1-h. 
 
\end{array}\right.
\end{equation} 
and $K_3$  defined by

\begin{equation}\left\{ \begin{array}{c}
 y_j=0, 1\leq j\leq n+1-h, \\
x_i y_j-x_j y_i=0, 1\leq i, j\leq n+1-h\\
x_i y_j-x_j y_i=0, n-h+2\leq i, j\leq n+1\\
(zx_0+x_{n+2-h}) y_j-(zy_0+y_{n+2-h}) x_j=0, 1\leq j\leq n+1-h,\\
x_j (zy_0+y_{n+2-h}) =0, n-h+2\leq j\leq n+1.
\end{array}\right.
\end{equation}

Pulling equations (3.6), (3.7) back to the variety $$\{\infty\}\times U\times X\times X,$$ and taking the closure, 
 we obtain 
 the cycle $ \Phi_1$  of the scheme 
$$\overline{ \tau^{-1} (K_2)\cap (\{\infty\}\times U\times X\times X )}.$$  
The other components from the scheme $\overline{ \tau^{-1} ( K_3)\cap (\{\infty\}\times U\times X\times X )}$
 lie in 
$$\{\infty\}\times \Upsilon\times X\times V^{n+1-h}.$$
At last we consider the components supported  in
$$\{\infty\}\times \{0\}\times X\times X.$$
So $\Omega_{\infty}^0$ is defined by 
\begin{equation}\left\{ \begin{array}{c}
x_i y_j=x_j y_i, 1\leq i, j\leq n+1-h\\
x_i y_j=x_j y_i, n-h+2\leq i, j\leq n+1\\
x_{n+2-h} y_j=y_{n+2-h} x_j, 1\leq j\leq n+1-h\\
 x_i y_j =0, 1\leq j\leq n+1-h, n-h+2\leq i\leq n+1, \\
x_j y_{n+2-h} =0, n-h+2\leq j\leq n+1.
\end{array}\right.
\end{equation}
Similarly we observe the fourth equations. There are two types of components.
The first one denoted by $F_1$
is defined by
$$
\left\{ \begin{array}{c}
x_{n+2-h}=\cdots=x_{n+1}=y_{n+2-h}=0 \\
x_i y_j=x_j y_i, 1\leq i, j\leq n+1-h\end{array}\right.
$$
and the second $F_2$ (for all $z$) lies in 
$$\mathbf P^{n+1}\times  \mathbf P^h$$
where $\mathbf P^h$ is the $h$-dimensional subspace defined by
$$x_i=0, 1\leq i\leq n+1-h.$$
Then the component of $\Theta_\infty$ 
originated from $F_2$ will be projected to $V^{n+1-h}$.
\par

The other component,  determined by $F_1$ has dimension $n-1$. So it can't be an irreducible component  of
$\Theta_\infty$ ( must be the boundary of the component  $ \Phi_1$). 
\bigskip

So overall, we proved that
\bigskip

\begin{proposition}
\begin{equation}
\Theta_\infty=\Phi_1+\Phi_2
\end{equation}
where $\Phi_2$ is a cycle lying in $$\{\infty\}\times \Upsilon\times X\times V^{n+1-h}.$$
Similarly, we'll denote the projection of $\Phi_i$ to $\Upsilon\times X\times X$ also by $\Phi_i$
\end{proposition}

\bigskip

\noindent $\bullet$  Case $t=0$.  \medskip

 Notice that $\Omega_t$ is defined by 

\begin{equation}\left\{ \begin{array}{c}
x_i y_j=x_j y_i, 1\leq i, j\leq n+1-h\\
x_i y_j=x_j y_i, n-h+2\leq i, j\leq n+1\\
y_i x_j=0, 1\leq j\leq n+1-h, n-h+2\leq i\leq n+1, \\
(zx_0+x_{n+2-h}) y_j=(zy_0+y_{n+2-h}) x_j, 1\leq j\leq n+1-h\\
y_j (zx_0+x_{n+2-h})=0, n-h+2\leq j\leq n+1.
\end{array}\right.
\end{equation}
Then the equations indicate the components of $\Theta_0$   are divided into two types: first type is defined by as the pullback of  
\begin{equation}
y_{n+1-h}=\cdots=y_{n+1}=0.
\end{equation}
The second type denoted by $W$ is   defined  
by the pullback of \begin{equation}
x_1=\cdots=x_{n+1-h}=0.
\end{equation}
 \bigskip

To summarize it

\begin{proposition}
\begin{equation} P_\ast (\Theta_0)=\zeta_1+\zeta_2,\end{equation}

where $\zeta_1$ lies in $W\times X$ and $\zeta_2$ lies in $X\times V^{h}$. \end{proposition}
\bigskip

\subsection{ Non-algebraic part of end cycles  }
\bigskip

In this subsection, we intersect the algebraic end cycles in the 
previous subsection with currents to obtain the end cycles in the cone family. 
We'll focus on the $\infty$-end cycle.
\bigskip

\noindent $\bullet$  The $1$-end cycle 
\bigskip

\begin{proof} of Proposition 1.6:  First we recall the proposition. 
 Let $\sigma$ be a singular cycle on $X$. Then we need to prove that with product De Rham data on $X\times X$ and projection De Rham data on $X$,  

\begin{equation}
\Psi_1 (\sigma)=m\sigma+\omega_\ast (\sigma)
\end{equation}
where $m$ is a natural number, where $(\cdot)_\ast$ denotes the correspondence of currents (see [6]).  

\bigskip
We notice that the factorization of the projection yields 
\begin{equation}
\Psi_1(\sigma)= ( \mathcal P_\ast (\Theta_1) )_\ast (\sigma).
\end{equation}

By Proposition 3.1
\begin{equation} ( \mathcal P_\ast (\Theta_1) )_\ast (\sigma)=m(\Delta_{X\times X})_\ast(\sigma)+\omega_\ast (\sigma)
\end{equation}

By Claim 2.11, [6], \begin{equation}
m(\Delta_{X\times X})_\ast(\sigma)=m\sigma.
\end{equation}
This completes the proof. 
\end{proof}

\bigskip

\noindent $\bullet$ $\infty$-end cycle
\bigskip

\begin{proposition} For any $\sigma\in \mathcal C(X)$,  
\begin{equation}
\Psi_{\infty} (\sigma)=(\Phi_1)_\star(\sigma)+(\Phi_2)_\star(\sigma)
\end{equation}
where $(\Phi_i)_\star$,  similar to (1.22),  is 
$$(\mu_4)_\ast [\Phi_i\wedge \{\infty\}\otimes \Upsilon\otimes \sigma\otimes X].$$
\end{proposition}

\bigskip

\begin{proof} This is the consequence of Proportion 3.2.

\end{proof}

\bigskip

In the following we calculate the currents $(\Phi_1)_\star(\sigma), (\Phi_2)_\star(\sigma)$. 
\bigskip

\begin{lemma}
\begin{equation}
(\Phi_2)_\star(\sigma)=0.\end{equation}
\bigskip

\end{lemma}
\bigskip

\begin{proof}

If $(\Phi_2)_\star(\sigma)\neq 0$, then it has dimension $\geq 2n-2h-1$. 
However $(\Phi_2)_\star(\sigma)$ is the inclusion  of a current in  $V^{n+1-h}$.  
Since $$dim(V^{n+1-h})=2n-2h $$ 
 $V^{n+1-h}$ can't have a non-zero current of dimension $> 2n-2h$.  
 We have a contradiction.
We complete the proof.

\end{proof}
\bigskip

Our observation to the Lefschetz standard conjecture is that  the factorization of the identity (1.6) 
will occur in an open set of the
projective variety, i.e in a quasi-projective variety for non-quotient geometric objects -- the $\infty$-end cycle.  For any algebraic cycle $A$, its intersection with the  open set 
is denoted by $\overset{\circ}{A}$.  Then we need a technique to pass on from a quasi-projective to a projective.  
This is the closure in the following definition. \bigskip

\begin{definition}
Let $X$ be a manifold and $U\subset X$ be an open set such that $\overline{ U}=X$. 
Let $T$ be a current on $U$. For any $\phi \in\mathscr D(X)$, if the evaluation
$\int_U \phi$ as an improper integral is convergent, then we define the current
$\bar T$ by the convergent evaluation
\begin{equation}
\int_{\bar T} \phi :=\int_U \phi.
\end{equation}
(The continuity of the functional also follows from formula (3.20) ) 

\end{definition}
\bigskip

\begin{ex} $X$ is a smooth projective variety over $\mathbb C$, and 
 $T$ is an irreducible subvariety. Let $U$ be a non empty Zariski open set
of $X$. Then currents' closure $\overline {T\cap U}=T$ provided  $T\cap U\neq \varnothing$.

\end{ex}
\bigskip

The following proposition shows that the passing on from quasi-projective to projective is cohomologically trivial. 
\bigskip

\begin{proposition}       For $\sigma\in \mathcal C_o(X)$,   
\par
(1) $$\overline{ [\overset{\circ}{\Phi_1}\wedge (U\otimes \sigma\otimes X)]}$$ \par\hspace{1CC} is well-defined  in 
 $\mathcal C_o(\Upsilon\times X\times X)$. \par

(2)  then the current 
\begin{equation}
 [\Phi_1\wedge (\Upsilon \otimes \sigma\otimes X)]
-\overline{ [\overset{\circ}{\Phi_1}\wedge (U\otimes \sigma\otimes X)]}
\end{equation}  \par\hspace{1CC} 
is  homologous to zero. 
\bigskip

Furthermore Proposition holds for arbitrary smooth projective variety $Y$ in the product 
$$\Upsilon\times Y\times X$$ with arbitrary algebraic cycle $\Phi_1$.
\end{proposition}

\bigskip

\begin{proof} Let $\Upsilon, X$ be equipped with De Rham data, and the subset $U\subset \Upsilon$ be equipped with the same  De Rham data as $\Upsilon$. Then $\Upsilon\times X\times X$ and $U\times X\times X$ are equipped with product De Rham data. Assume $X$, the second copy of $X$  in  $U\times X\times X$ is changed to the
 projection De Rham data with respect to the product De  Rham data of $U\times X\times X$. That projection De Rham data must be the same 
projection De Rham data with respect to the product De Rham data on  
$\Upsilon\times X\times X$.  In particular,  projection formulas in Proposition 2.8, [6]  for currents' projections 
 $$P_2, \Upsilon\times X\times X\to X(2nd\ factor),  P_2: U\times X\times X\to X(2nd\ factor)$$ hold.

\par Part (1). 
Let $\phi\in \mathscr D(\Upsilon\times X\times X)$. Let $\overset{\circ}{\phi}$ be the distribution whose restriction  to $U\times X\times X$ is $\phi$ and $0$ 
at $\{\propto\}\times X\times X$. 
Then there is a family of currents $\phi_n, n\in\mathbb N$,   such that  $\phi_n\in  \mathscr D(\Upsilon\times X\times X)$ are
$C^\infty$, supported in $U\times X\times X$ 
 and $\displaystyle{\lim_{n\to \infty}}\phi_n= \overset{\circ}{\phi}$
in  $\mathscr D'(\Upsilon\times X\times X)$. 
Let $P_2: U\times X\times X\to X(2nd \ factor)$ be the projection. By the definition, 
$\overline{ [\overset{\circ}{\Phi_1}\wedge (U\otimes \sigma\otimes X)]}$ evaluated at $\phi$ is defined to be the limit 
$$ \displaystyle{\lim_{n\to \infty}}\int_{[\overset{\circ}{\Phi_1}\wedge (U\otimes \sigma\otimes X)]}\phi_n.$$ So  
We calculate
\begin{align}
 & \displaystyle{\lim_{n\to \infty}}\int_{[\overset{\circ}{\Phi_1}\wedge (U\otimes \sigma\otimes X)]}\phi_n\\
&=\displaystyle{\lim_{n\to \infty}}\displaystyle{\lim_{\epsilon\to 0}}\int_{\overset{\circ}{\Phi_1}}
(P_2)^\ast (R_\epsilon(\sigma))\wedge \phi_n\\
& 
(\text {Change the improper integral to proper integral.})\nonumber\\
&
=\displaystyle{\lim_{n\to \infty}}\displaystyle{\lim_{\epsilon\to 0}}\int_{\Phi_1}
  (P_2)^\ast (R_\epsilon(\sigma))\wedge \phi_n
\\
&=\displaystyle{\lim_{n\to \infty}} \int _{(P_2)_\ast [ {\Phi_1}\wedge \phi_n]\wedge \sigma} 1.
\end{align}
Next we calculate the family of currents ${(P_2)_\ast [{\Phi_1}\wedge \phi_n]\wedge \sigma}$.
We should note the limit of currents $(P_2)_\ast ({\Phi_1}\wedge \phi_n)$ exists  as the functional 
$$f\to \int_{{\Phi_1}}\phi\wedge (P_2)^\ast f,$$for any test form  $f\in \mathscr D(X)$. 
This functional is equal to
$$(P_2)_\ast \overline {[\overset{\circ} {\Phi_1}\wedge \overset {\circ}{\phi}]}.$$
The following is the argument of the assertion.   
By the same result of Lemma 4.1 in [6], 
$(P_2)_\ast \overline {[\overset{\circ} {\Phi_1}\wedge \overset {\circ}{\phi}]}$ is Lebesgue.  Then
\begin{align}
 & \int_{{(P_2)_\ast ( [{\Phi_1}\wedge \phi_n]}-
(P_2)_\ast \overline {[\overset{\circ} {\Phi_1}\wedge \overset {\circ}{\phi}]}} f
=\int_{\Phi_1} (\phi_n-\phi)\wedge (P_2)^\ast (f).
\end{align}
So  (3.26) as $n\to \infty$ converges to $0$ uniformly on bounded set of forms $f$. Hence 
 ${(P_2)_\ast  [{\Phi_1}\wedge \phi_n]}$ converges to 
 $(P_2)_\ast \overline {[\overset{\circ} {\Phi_1}\wedge \overset {\circ}{\phi}]}$ in the space of currents. 
By the particular type of continuity
\footnote{There is no general claim for the continuity of intersection.  For instance it requires flatness in case of algebraic geometry.} of the intersection -- Proposition 4.12, [6], 
${(P_2)_\ast [{\Phi_1}\wedge \phi_n]}\wedge \sigma$ converges to 
$(P_2)_\ast \overline {[\overset{\circ} {\Phi_1}\wedge \overset {\circ}{\phi}]}\wedge \sigma$. Hence the number (3.25) converges to
\begin{equation}
 \int _{(P_2)_\ast \overline {[\overset{\circ} {\Phi_1}\wedge \overset {\circ}{\phi}]}\wedge \sigma} 1.
\end{equation}

\par
Next we show it is  closed. We have  
\begin{equation}
d\overline{ [\overset{\circ}{\Phi_1}\wedge (U\otimes \sigma\otimes X)]}
=\overline{ [ d\overset{\circ}{\Phi_1}\wedge (U\otimes \sigma\otimes X)]}\pm 
\overline{ [ \overset{\circ}{\Phi_1}\wedge (dU\otimes \sigma\otimes X)]}
.\end{equation}
\par Since $\overset{\circ}{\Phi_1}, U$ both are quasi-projective, the currents of integration over them are closed. So 
$$d\overline{ [\overset{\circ}{\Phi_1}\wedge (U\otimes \sigma\otimes X)]}=0.$$
This proves part (1). \par

Part (2).   Since $\overline{ [\overset{\circ}{\Phi_1}\wedge (U\otimes \sigma\otimes X)]}$ is well-defined,
 it suffices to evaluate the current of (3.21)  at
the closed forms.  
By K\"unneth decomposition it suffices to valuate at two $C^\infty$ forms, 
$$1\otimes \omega, \psi \otimes \omega$$
where $\omega$ is a closed form on $X\times X$ and $\psi$ is a $C^\infty$ 2-form supported in a neighborhood of a generic point in $\Upsilon\simeq \mathbf P^1$ and Poicar\'e dual to the point. 
Then
\begin{align}
 &\int_{\overline{ [\overset{\circ}{\Phi_1}\wedge (U\otimes \sigma\otimes X)]}}\psi\otimes \omega\\
&= \int_{ [\overset{\circ}{\Phi_1}\wedge (U\otimes \sigma\otimes X)]}\psi\otimes \omega\\
&= \int_{[\Phi_1\wedge (\Upsilon\otimes \sigma\otimes X)]}\psi\otimes \omega.\end{align}
For the form $1\otimes \omega$, we calculate
\begin{align}
   &\int_{\Phi_1\wedge (\Upsilon\otimes\sigma \otimes X)}1\otimes \omega\nonumber\\
 &=\displaystyle{\lim_{\epsilon\to 0}}\int_{\overset{\circ}{\Phi_1}} (P_2)^\ast( R_\epsilon(\sigma))\wedge (1\otimes \omega)\\
& (\text{ (3.32) only holds for the particular type  of  forms in $1\otimes \omega$.})\\ 
& =\int_{\overline { [\overset{\circ}{\Phi_1} \wedge (U\otimes \sigma\otimes X)]}} 
1\otimes \omega
\end{align}
where (3.34) is an improper integral. But in part (1), we have shown 
 \begin{equation}
\displaystyle{\lim_{n\to \infty}}\int_{[\overset{\circ}{\Phi_1}\wedge (U\otimes \sigma\otimes X)]}(\bullet) 
\end{equation}
 is convergent. 
 So  we obtain that 
\begin{equation}
\int_{ [\Phi_1\wedge (\Upsilon \otimes \sigma\otimes X)]
-\overline{ [\overset{\circ}{\Phi_1}\wedge (U\otimes \sigma\otimes X)]}}\alpha=0\end{equation}
for any closed form $\alpha$ on $\Upsilon\times X\times X$. So we complete the proof of part (2). 
Furthermore the argument still holds when the second factor $X$ in $\Upsilon\times X\times X$ is replaced by  arbitrary 
smooth projective variety $Y$ and  algebraic cycle $\Phi_1$. 

\end{proof}

\bigskip

\begin{proof} of Proposition 1.7: 
By Proposition 3.2, 
$\Phi_1$ is the projection of the closure of
the zero locus 
\begin{equation}
\Gamma
\end{equation}
on quasi-projective variety
\begin{equation}
 {\infty}\times U\times X\times X
\end{equation}
of restricted sections

\begin{equation}\left\{\begin{array}{c}
\alpha_{ij}=\tau^\ast (x_iy_j-x_jy_i), 1\leq i, j\leq n+1-h, \\
\alpha_{ij}=\tau^\ast (x_iy_j-x_jy_i), n+2-h\leq i, j \leq n+1\\
\alpha_j=\tau^\ast \biggl((zx_0+x_{n+2-h}) y_j-(zy_0+y_{n+2-h}) x_j\biggr), 1\leq j\leq n+1-h\\
\tau^\ast (x_i), n+2-h\leq i\leq n+1.
\end{array}\right. \end{equation}
where $\tau^\ast (\bullet )$ are regarded as sections of the line bundle over ${\infty}\times U\times X\times X$.

The following  is the key observation of equations (3.39) for $\Gamma$: 
$\Gamma$  is a local complete intersection of codimension $n+1+p$. Let's see this. 
We divide the collection of equations (3.39) into two parts.
First let
$$
\overset{\circ} I_h\subset \{\infty\}\times U\times X\times X$$ be the zero locus of
\begin{equation}\left\{\begin{array}{c}
\alpha_{ij},  1\leq i, j\leq n+1-h, \\
\alpha_{ij}, n+2-h\leq i, j \leq n+1\\
\alpha_j, 1\leq j\leq n+1-h.\\
\end{array}\right. \end{equation}
This portion defines an l.c.i. on an open set of codimension $n+1$. 
Notice that  $\overset{\circ} I_h$ has been  defined in introduction but through a different expression.  Secondly
we also notice 
$$\{\infty\}\times U\times V^h\times X$$ is the zero locus 
of 
\begin{equation}
\tau^\ast (x_i), n+2-h\leq i\leq n+1
\end{equation}
where $x_i$ are the $\tau^\ast(\mathcal O_{\mathbf P^{n+1}}(1))$ sections over the
first copy of $\mathbf P^{n+1}$.  This portion defines a complete intersection of codimension $h$. 
Hence the grouping of the equations (3.39) determines  that the cycle $ \Phi_1$ is the cycle of the scheme
\begin{equation}
\overline {\overset{\circ} I_h\cap (\{\infty\}\times U\times V^h\times X)\cap (\{\infty\}\times U\times X\times X)}.
\end{equation}
Now we observe the intersection $(\{\infty\}\times U\times V^h\times X)\cap (\{\infty\}\times U\times X\times X)$ is proper
in $\{\infty\}\times U\times X\times X$, and 
the intersection $$\overset{\circ} I_h\cap \biggl( (\{\infty\}\times U\times V^h\times X)\cap (\{\infty\}\times U\times X\times X)\biggr)$$ is also proper 
in $\{\infty\}\times U\times V^h\times X$. Then we change the scheme intersection to cycle intersection to have the triple intersection 
\begin{equation}
\Phi_1=
\overline{\overset{\circ} I_h\cdot \biggl ((\{\infty\}\times U\times V^h\times X)\cdot  (\{\infty\}\times U\times X\times X)\biggr) }
\end{equation}
where the two static intersections occur in two different quasi-projective varieties: $\{\infty\}\times U\times X\times X$ and $\{\infty\}\times U\times V^h\times X$.
Applying the associativity in the intersection of currents, we obtain that the equality from the quasi-projective variety
$${\infty}\times U\times X\times X, $$

\begin{equation}
(\Phi_1)_\star (\sigma)=
\overline{\overset{\circ} {I_h}\cdot \biggl ((\{\infty\}\otimes U\otimes [\sigma\wedge V^h]\otimes X)\biggr)}\\
\end{equation}
Therefore by Proposition 3.8, there is homologically trivial cycle $\mathcal T_\sigma$ such that 
\begin{equation}\begin{array}{c}
\Psi_\infty(\sigma)=(\eta_4)_\ast \biggl (
\overline {\overset{\circ} I_h\cdot (\{\infty\}\times U\times (\sigma\cdot V^h)\times X)}\biggr)+\mathcal T_\sigma\\
=(m+d) Con_h\circ v^h(\sigma)+\mathcal T_\sigma. \end{array}
\end{equation}

We complete the proof of Proposition 1.7. 

\end{proof}

\bigskip

\noindent $\bullet$ The $0$-end cycle 

\bigskip

\begin{proof} of Proposition 1.8:  First we 
know 
$$P_\ast (\Theta_0)=\zeta_1+\zeta_2.$$

Let $\sigma\in \mathcal C_o(X)$. If $dim(\sigma)<2n-2h$, then  
$W\cap supp(\sigma)=\varnothing$.  Therefore $(\zeta_1)_\ast(\sigma)=0$.
Applying Proposition 4.5, [6] concerning the support of intersection, 
we obtain that
$$supp((\zeta_2)_\ast(\sigma))\subset V^h.$$
we complete the proof.

\end{proof}

\bigskip

\bigskip

\begin{center}\section{Algebraicity and cohomologicity }\end{center}

\subsection{ Algebraicity}

\bigskip

\begin{definition} Let $X$ be a smooth projective variety of dimension $n$ over $\mathbb C$.
Let $T\in (\mathcal D')^{2p+k} (X)$ be a real current of degree $2p+k$.  If the support of $T$ lies
in  an algebraic set of codimension $p$, we say $T$ has an algebraic level $k$.

\end{definition}
\bigskip

{\bf Remark} The algebraic level of a current $T$ is not a unique number. But there is the minimum. 
\bigskip

\begin{proposition} Let $\sigma\in \mathcal C_o(V^h)$. 
If $\sigma$ has an algebraic  level $k$, 
so does  $Con_h(\sigma)$. 
 
\end{proposition}

\bigskip

{\bf Remark} This is the algebraicity of $Con_h$.  

\bigskip

\begin{proof} 
Let $\sigma\in \mathcal C_o^{2n-i}$ of dimension $i$. 
As a current, \begin{equation}
Con_h(\sigma)={1\over m+d}(\eta_4)_\ast\overline  {\biggl[\overset{\circ}{I_h}
\wedge (\{\infty\}\times U\times \gamma \times X)\biggr ] }.
\end{equation} has
dimension $i+2h$.   So if $\sigma$ has the algebaric level $k$, 
 we need to find an algebraic set of dimension  $h+{k+i\over 2}$ containing $Con_h(\sigma)$. 
By the assumption,
 there is an algebraic cycle $B$ of dimension  ${i+k\over 2}$ in $X$  such that 
$$supp(\sigma)\subset B.$$
If $\Upsilon$ is generic with respect to $B$, then the algebraic cycle
\begin{equation}
I_h\cdot (\{\infty\}\times U\times B \times X)\end{equation}
is well defined,   of dimension $h+{k+i\over 2}$. Applying Proposition 4.5, [6] concerning the support of intersection, we obtain that
\begin{equation}
(\eta_4)_\ast (I_h\cdot (\{\infty\}\times U\times B \times X))\end{equation}
is an $h+{k+i\over 2}$  dimensional algebraic cycle 
containing  the support of  $Con_h(\sigma)$.
Since the cohomology has a finite dimension,  generic and  fixed curve $\Upsilon$ 
can be selected  for all $Con_h(\sigma)$ to have level $k$.

\par

\end{proof}

\bigskip

\subsection{Cohomology with support}

\bigskip

\noindent $\bullet$\quad {\bf 1-end cycle}
\bigskip

\begin{proposition}
For any natural number $i\neq n$, the currents' correspondence 
  $ \omega_\ast$ descends to  the isomorphism on the cohomology 
\begin{equation}\begin{array}{ccc}
H^i(X;\mathbb Q) &\rightarrow & H^i(X;\mathbb Q)\\
\theta &\rightarrow & d\theta
\end{array}\end{equation}
where $d=deg(\mu^\ast(\mathcal O_{\mathbf P^{n+1}}(1)))$. 

\end{proposition}
\bigskip

\begin{proof} Applying the cohomologicity of currents' intersection, Property 2.6, [6], 
For a singular cycle $\sigma$ of codimension $i$, 
\begin{equation}
\langle \omega_\ast (\sigma)\rangle =\langle \omega\rangle_\ast \langle \sigma\rangle.
\end{equation}
Hence it suffices to prove Proposition 4.3 for the cohomological correspondence $\langle \omega\rangle_\ast$. 
First we consider the case $i<n$. 
 Let $\mathbf P^{n}$ be a  hyperplane of
$\mathbf P^{n+1}$ defined by $x_0=0$.   Let $\theta\in H^{i}(X;\mathbb Q)$.
$\langle \omega\rangle _\ast(\theta)$ is well-defined cycle class in $H^{i}(X;\mathbb Q)$.
Let $\phi\in H^{2n-2-i}(X;\mathbb Q)$. Let's
calculate $$\biggl( \langle \omega\rangle _\ast(\theta)\cup \langle V^1\rangle, \phi\biggr) $$
where $V^1$ is the hyperplane section $div(\mu^\ast(x_0))$, and $(\bullet, \bullet)$ is the intersection number. 
We use De Rham's theorem to compute the intersection in the cohomology of real coefficients. So we  assume $\theta, \phi$ are  closed smooth forms on $X$.
Let \begin{equation}\xymatrix{
 &X\times X \ar[dl]_{Pr_1}\ar[dr]^{Pr_2} \\
X & &X 
}\end{equation}
be the two projections. 
Then by De Rham's theorem  \begin{equation}
\biggl( \langle \omega\rangle _\ast(\theta)\cup \langle V^1\rangle, \phi\biggr)=
\int_{\omega\cdot (X\times V^1)} Pr_1^\ast (\theta)\wedge Pr_2^\ast (\phi).\end{equation}
The right hand side of (4.7) is the integral over the algebraic cycle $\omega\cdot (X\times V^1)$. 
There is covering map of degree $d$  in $\mathbf P^{n+1}$
\begin{equation}
\begin{array}{ccc}
\mathcal K: X &\rightarrow & \mathbf P^n.\end{array}
\end{equation} By the general position assumption on $X$,  the restriction 

\begin{equation}
\begin{array}{ccc}
\mathcal K|_{V^{-1}}: V^{-1} &\rightarrow & V^1\end{array}
\end{equation}
is also a covering map of degree $d$, where $V^{-1}\subset \mathcal K^{-1} (V^{1})$ is the component dominating $V^1$. 
We found that
$$\omega\cdot (X\times V^1)$$ is the cycle of the  scheme
$$ \omega\cap (X\times V^1).$$ Thus the covering map $\mathcal K$
 gives the isomorphic diagram
\begin{equation}\xymatrix{
\omega\cap (X\times V^1)\ar[dr]_{\rho} &\simeq  &V^{-1}\ar[dl]^{\mathcal K}\\
& V^1  }\end{equation}
where $\omega\cap (X\times V^1)$ is isomorphic to the graph of $\mathcal K|_{ V^{-1}}$, and $\rho$ is the projection
from the graph of the map to the image of the map. 
Then we can use the computation in currents' evaluation where algebraic varieties 
$\omega, V^1$ and the differential form $\theta$ are all currents. 
\begin{equation}\begin{array}{c}
\int_{\omega\cdot (X\times V^1)} Pr_1^\ast (\theta)\wedge Pr_2^\ast (\phi)=
d \int_{V^1} \theta\wedge \phi=d\int_{V^1\wedge \theta}\phi
\end{array}\end{equation}
where $V^1\wedge \theta$ is the intersection of currents.
The computation (4.11) shows a cohomological equality
\begin{equation}
\langle \omega\rangle _\ast(\theta)\cup \langle V^1\rangle=k\langle\theta\rangle \cup \langle  V^1 \rangle.
\end{equation}
 Therefore 
\begin{equation} L (\langle \omega\rangle _\ast(\theta))= L (k \theta),
\end{equation}
where $L$ is the Lefschetz operator defined in (1.1). 
Applying the hard Lefschetz theorem, we obtain that  
\begin{equation}
\langle \omega\rangle_\ast  ( \theta)=d \theta.\end{equation}

\par
Next we assume $i>n$. Because $\omega$ is symmetric, 
the transpose $\langle \omega\rangle^\intercal_\ast $ is equal to $\langle \omega\rangle_\ast$.
Now we let $\theta, \alpha$ be cohomological classes of degrees $i$ and $2n-i$.
 We calculate the intersection number, 
\begin{equation}\begin{array}{c}
(\langle \omega\rangle_\ast  \theta, \alpha)=
(\theta, \langle \omega\rangle^{\intercal}_\ast \alpha)\\
=(\theta, \langle \omega\rangle_\ast  \alpha)\\
=( \theta,  d\alpha)\\
=(d\theta, \alpha).\end{array}\end{equation}
Therefore 
\begin{equation}
\biggl(\langle \omega\rangle_\ast \theta-d \theta,  \alpha\biggr)=0
\end{equation}
holds for all $\theta, \alpha$. 
We obtain that  
$$\langle \omega\rangle_\ast \theta=d \theta.$$
 We complete the proof for all cases. 
\end{proof}

\bigskip

\noindent  $\bullet$\quad {\bf $\infty$-end cycle}
\bigskip

\begin{lemma} Assume $V^h\times X$ is equipped with a product De Rham data. 
When restricted to the closed Lebesgue currents $$\mathcal C_o^{n+h} (V^h)$$
the operator $Con_{h, p}$ is cohomological for the homological equivalence of $X$
, i.e. it sends an $n-h$ dimensional current of $V^h$,  exact in $X$  to an $n+h$ dimensional current of $X$,  exact in $X$.
\end{lemma}

\bigskip

\begin{proof} Let $\lambda\in (\mathcal C_o)(X)$ be the pushforward of current  
of $\lambda_h\in \mathcal C_o (V^h)$ under the inclusion 
map.
  \bigskip

\begin{ass}
\begin{equation}
  v^h\circ Con_{h, p}(  \lambda_h)=l\lambda_h+\mathcal J_{\lambda_h},
\end{equation}\end{ass}
where $l$ is a natural number and  $\mathcal J_{\lambda_h}$ is homologous to zero.

Let  $Pr_3: \Upsilon\times V^h\times X\to X$ be the projection.  Let  the Cartesian product $\Upsilon\times X\times X$
have a product De Rham data and 
and the last factor $X$ have projection De Rham data.  
By the projection formula ( Proposition 2.8, [6]), 

\begin{align}
   & V^h\wedge (Pr_3)_\ast \overline { [\overset{\circ}{I_h}\wedge (U\otimes \lambda_h\otimes X)]}\\
 & =(Pr_3)_\ast \biggl [ (\Upsilon\otimes V^h\otimes V^h)\wedge \overline { [\overset{\circ}{I_h}
\wedge (U\otimes \lambda_h\otimes X)]}\biggr].
\end{align}

By Proposition 3.8, which is also valid in this case, 
$$\overline { [\overset{\circ}{I_h}\wedge (U\otimes \lambda_h\otimes X)]}$$ is homologous to 
$$ [{I_h}\wedge (\Upsilon\otimes \lambda_h\otimes X)].$$
Hence (4.19) is equal to

\begin{equation}
(Pr_3)_\ast \biggl[ ( [ \Upsilon\otimes V^h\otimes V^h)\wedge [{I_h}\wedge (\Upsilon\otimes \lambda_h\otimes X)]\biggr] +\mathcal J_{\lambda_h},
\end{equation}
where $\mathcal J_{\lambda_h}$ is homologous to zero.
Then the associativity says
\begin{align}
 &(Pr_3)_\ast \biggl[ [ \Upsilon\otimes V^h\otimes V^h)\wedge [{I_h}\wedge (\Upsilon\otimes \lambda_h\otimes X)]\biggr] \\
& = (Pr_3)_\ast\biggl[ ( \Upsilon\otimes V^h\otimes V^h)\wedge {I_h}\wedge (\Upsilon\otimes \lambda_h\otimes X)]\biggr]
\end{align}
We first consider the current $( \Upsilon\otimes V^h\otimes V^h)\wedge {I_h}$ which is
the current of integration over the algebraic cycle
$$( \Upsilon\times V^h\times V^h)\cdot {I_h}.$$
By the definition of the variety $I_h$,  the projection  
\begin{equation}\begin{array}{ccc}
Pr_{23}: ( \Upsilon\times V^h\times V^h)\cdot {I_h}  &\rightarrow & V_h\times V_h
\end{array}\end{equation}
is a multiple cover of the diagonal $\Delta_{V^h}$ whose multiplicity is denoted by $l$.
Let $$Pr: V^h\times V^h\to V^h(2nd\ factor)$$ be the projection. 
Using the calculation of  the intersection of currents, specifically,  projection formula in Proposition 2.8, and claim 2.11 in [6], 
we obtain that 
\begin{align}
 &(Pr_3)_\ast \biggl[( [ \Upsilon\otimes V^h\otimes V^h)\wedge [{I_h}\wedge (\Upsilon\otimes \lambda_h\otimes X)]\biggr]\\
&=(Pr)_\ast  [l\Delta_{V^h}\wedge( \lambda_h\otimes X)]\\
&= l \lambda_h.\end{align}
Combining with (4.20), we complete the proof of Claim 4.5. 
Now we assume  $\lambda$ has degree $n+h$, i.e. restrict $Con_h$ to $n-h$ dimensional  cycles.   Then 
we can  apply hard Lefschetz theorem to obtain that 
if $\lambda$ is cohomologically trivial in $X$, so is  $ Con_{h, p}( \lambda_h)$ in
$X$. This shows that $Con_{h, p}$ is cohomological. We complete the proof. 
\end{proof}

\bigskip

\begin{proof} of Proposition 1.11: Lemma 4.4 shows that $Con_{h, p}$ is cohomological (with respect to the homological equivalence of $X$), i.e. 
$Coh_{h, p}$ is reduced to a homomorphism $\langle Coh_{h, p}\rangle$ on the homogeneous part of cohomology,
\begin{equation}\begin{array}{ccc}
\langle Coh_{h, p}\rangle: N^pH^{2p+k}_{V_h}(X) &\rightarrow & N^qH^{2q+k}(X)
\end{array}\end{equation}
where $N^pH^{2p+k}_{V_h}(X)$ is the cohomology supported on $V^h$.  Furthermore (1.29) becomes
\begin{equation}
\langle Coh_{h, p}\rangle\circ L_k^h=identity.
\end{equation}

Also by Claim 4.5  for $\delta\in \mathcal C_o(X)$,  we have the cohomological version 

\begin{equation}
u^h\cup  \langle Con_{h, p} (\delta)\rangle = l_\delta \langle \delta\rangle .
\end{equation}
for some rational number $l_\delta$.   In the following we prove that $l_\delta=1$. 
Let $$\sigma\in \mathcal C_o(X).$$   We have 
\begin{equation}
u^h\cup \langle Con_{h, p}  \circ  v^h ( \sigma )\rangle= l_{v^h(\sigma)}  u^h\cup \langle \sigma\rangle
\end{equation}
and on the other hand by (4.28) 
\begin{equation}
u^h \cup \langle Con_{h, p}\circ  v^h(\sigma )\rangle =\langle v^h(\sigma )\rangle , \quad in\ X
\end{equation}
 We obtain that 
$$ l_{v^h(\sigma)}  u^h\cup \langle \sigma\rangle = \langle v^h(\sigma )\rangle. $$
By the hard Lefschetz theorem for $L^h$, $l_{v^h(\sigma)}=1$, and
$l_\delta=1$ holds for all closed $\delta$.  
This result has two implications: \par
  \begin{equation}\begin{array}{c}\langle Con_{h, p}\rangle\circ L_k^h=identity, \ by\ (4.28)\\
L_k^h \circ \langle Con_{h, p}\rangle=identity,  \ by\ (4.29).\end{array}\end{equation}
Hence $\langle Coh_{h, p}\rangle$ is an isomorphism 

\begin{equation}
 N^pH^{2p+k}_{V_h}(X) \simeq  N^qH^{2q+k}(X).\end{equation}

\end{proof}
\bigskip

\subsection{Cohomology without support}

\bigskip

\noindent  {\bf $0$-end cycle}\bigskip

\begin{proof} of Proposition 1.12: At $t=0$, we consider $\sigma$ having dimension $2p+k>0$. 
By Proposition 1.8, 
\begin{equation}
\Psi_0 (\sigma),
\end{equation}
must lie in $V^h$.  By (1.35), a rational multiple of $\sigma$  is homologous  to 
\begin{equation}
\Psi_0 (\sigma),
\end{equation}
which  a cycle of the same dimension,  lying in $V^h$. To consider the level, 
 we let $A$ be an algebraic cycle containing $\sigma$.  The same operation 
$\Psi_0(\cdot)$ takes $A$ to an algebraic cycle $\Psi_0(A)$ of the same dimension, lying in $V^h$. 
Therefore the level is not changed. We have 

Therefore \begin{equation}
N^q H^{2q+k}_{V^h}(X)\simeq N^q H^{2q+k}(X).
\end{equation}
 We complete the proof of Main theorem.

\end{proof}

\bigskip

\begin{appendices}

\section{Coniveau Filtration of currents' version}\par
While we review 
the definition of coniveau filtration, we'll give another 
description using currents.  Recall that in [4], Grothendieck created a  filtration $Filt'{^p}$, 
which he called ``Arithmetic filtration, as it embodies
deep arithmetic properties of the scheme ". 
   This later was referred to as the coniveau filtration, denoted by
$$N^pH^{2p+k}(X)$$ 
where $p$ is  the coniveau and $k$ is  the level.
It is defined as a linear span of kernels of the linear maps
\begin{equation}\begin{array}{ccc}
H^{2p+k}(X;\mathbb Q) &\rightarrow &H^{2p+k}(X- W;\mathbb Q)
\end{array} \end{equation}
for all subvarieties  $W$ of codimension at least $p$. 
This is the original definition. 
In the same paper, Grothendieck  immediately interpreted it 
as a linear span of images of Gysin homomorphisms
\begin{equation}\begin{array}{ccc}
H^{dim(W)+2p+k-2n}(\tilde W;\mathbb Q) &\rightarrow &H^{2p+k}(X;\mathbb Q)
\end{array} \end{equation}
for all subvarieties  $W$ of codimension at least $p$ with the smooth resolution $\tilde W$.  
 Now it is known that the proof of the description requires Deligne's  mixed Hodge  structures. 
  We introduce another interpretation. It is  through currents, which are known to unite
both homology and cohomology.  Let $\mathcal D'(X)$ be the space of currents over $\mathbb R$ on $X$. 
 Let $C\mathcal D'(X)$ be its subset of closed currents and $E\mathcal D'(X)$ be its subset of exact currents.
Then 
\begin{equation} {C\mathcal D'(X)\over   E\mathcal D'(X)}= \sum H^{\bullet} (X;\mathbb R).\end{equation}
There is a restriction map on currents
\begin{equation}\begin{array}{ccc}
\mathcal R: \mathcal D'(X) &\rightarrow \mathcal D'(X-W)
\end{array}\end{equation}
for a subvariety $W$.

Using (A.3) and (A.4),  we define  $$\mathcal D^{p} H^{2p+k}(X)$$
to be the linear span of classes in $H^{2p+k} (X;\mathbb Q)$ such that they lie in
\begin{equation}
{C\mathcal D'(X)\cap kernel(\mathcal R)\over   E\mathcal D'(X)\cap kernel(\mathcal R)}.
\end{equation}
 for some $W$ of codimension at least $p$.
 \bigskip

\begin{proposition}
Let $ X$ be a smooth projective variety over $\mathbb C$. 
Then 
\begin{equation}
\mathcal D^p H^{2p+k}(X)=N^p H^{2p+k}(X).
\end{equation} It says that 
the cohomology class $\alpha$ lies in
\begin{equation}
 N^p H^{2p+k}(X) 
\end{equation}
if and only if it is represented by a current whose support is contained in an algebraic set of codimension at least $p$.

\end{proposition}

\begin{proof}: See [8].
\end{proof}

\bigskip

\end{appendices}

\end{document}